\documentclass[11pt]{article}
\usepackage{amsmath}
\usepackage{amsfonts}
\usepackage{amsthm}
\usepackage{verbatim}

\begin{document}

\def\N{\mathbb{N}}
\def\F{\mathbb{F}}
\def\Z{\mathbb{Z}}
\def\R{\mathbb{R}}
\def\Q{\mathbb{Q}}
\def\H{\mathcal{H}}

\parindent= 3.em \parskip=5pt

\centerline{\bf{On particular families}} 
\centerline{\bf{of hyperquadratic continued fractions}}
\centerline{\bf{in power series fields of odd characteristic}}
\vskip 0.5 cm
\centerline{\bf{by A. Lasjaunias}}
\vskip 0.5 cm
\noindent {\bf{Abstract.}} We discuss the form of certain algebraic continued fractions in the field of power series over $\F_p$, where $p$ is an odd prime number. This leads to give explicit continued fractions in these fields, satisfying an explicit algebraic equation of arbitrary degree $d\geq 2$ and having an irrationality measure equal to $d$. Our results are based on a mysterious finite sequence of rational numbers.  
\vskip 1 cm

\noindent{\bf{1. Introduction}} 
\par We are concerned with continued fractions in the fields of formal power series over a finite field. For a general account on this matter the reader may consult Schmidt's article \cite{S} and also Thakhur's book \cite{T}. Let $p$ be a prime number and $\F_q$ the finite field of characteristic $p$, having $q$ elements. Given a formal indeterminate $T$, we consider the ring $\F_q[T]$, the field $\F_q(T)$ and $\F_q((T^{-1}))$, here simply denoted by $\F(q)$, the field of power series in $1/T$ over the finite field $\F_q$. A non-zero element of $\F(q)$ can be written as
 $$\alpha=\sum_{i\leq i_0}u_iT^i \quad \text{ where }\quad i\in \Z,\quad u_i\in \F_q \quad \text{ and }\quad u_{i_0}\neq 0.$$
An ultrametric absolute value is defined over this field by $\vert 0 \vert =0$ and $\vert \alpha \vert =\vert T \vert^{i_0}$ where $\vert T\vert$ is a fixed real number greater than 1. We also consider the subset $\F(q)^+=\lbrace \alpha \in \F(q) \quad s.t. \vert \alpha \vert >1\rbrace$. Note that $\F(q)$ is the completion of the field $\F_q(T)$ for this absolute value. 
\par In power series fields over a general finite field $\F_q$, where $q$ is a power of $p$, contrarily to the case of real numbers, the continued fraction expansion for many algebraic elements can be explicitly given. This phenomenon is due to the existence of the Frobenius isomorphism in these fields. A particular subset of $\F(q)$, denoted by $\H(q)$, containing certain algebraic elements called hyperquadratic, has been considered (see \cite{BL}). Let $t\geq 0$ be an integer and $r=p^t$, an irrational element of $\F(q)$ is called hyperquadratic of order $t$ if it satisfies a non-trivial algebraic equation of the following form : $$uX^{r+1}+vX^r+wX+z=0 \quad \text{ where }\quad (u,v,w,z)\in (\F_q[T])^4.$$ Note that a hyperquadratic element of order $0$ is simply irrational quadratic. Irrational quadratic power series, as quadratic real numbers, have an ultimately periodic continued fraction expansion. However $\H(q)$ contains power series of arbitrary large algebraic degree over $\F_q(T)$ and the continued fraction expansion for various elements in this class has also been given explicitly. The consideration of this subset was first put forward in the study of Diophantine approximation, beginning with Mahler's article \cite{M}. For a survey on this topic and more referecences, the reader may also consult \cite{L}. 
\par Let us recall that, for an irrational element $\alpha \in \F(q)$, the irrationality measure is defined by :
$$\nu(\alpha)=-\limsup_{\vert Q \vert \to \infty}(\log \vert \alpha -P/Q \vert/\log \vert Q \vert),$$
where $P$ and $Q$ belong to $\F_q(T)$. Then $\nu(\alpha)$ is a real number greater or equal to $2$. By adapting a theorem on rational approximation for real numbers, due to Liouville in the 19-th century, Mahler \cite{M} proved that, if $\alpha$ is an algebraic element of degree $d>1$ over $\F_q(T)$, then we have $\nu(\alpha)\in [2;d]$. Furthermore if $\alpha$ is any irrational number in $\F(q)$, having the infinite continued fraction expansion $\alpha=[a_1,a_2,\dots,a_n,\dots]$, then the irrationality measure is directly connected to the sequence of the degrees of the partial quotients (see \cite[p.~214]{L}) and we have
$$\nu(\alpha)=2+\limsup_{n\geq 1}(\deg(a_{n+1})/\sum_{1\leq i\leq n}\deg(a_i)).$$
\par In this note, we present certain particular algebraic continued fractions in $\F(p)$ with odd $p$, which can be fully described. In section 2, we recall several generalities on continued fractions. In section 3, we describe a large family of algebraic continued fractions and, in section 4, we exhibit two particular sub-families where the sequence of partial quotients is regularly distributed. The importance of these last continued fractions is highlighted in the last section, bringing families of hyperquadratic elements, in $\F(p)$ with odd $p$, having a prescribed algebraic degree, an explicit continued fraction expansion and an irrationality measure equal to the algebraic degree. 
\par The reader will observe that this paper contains mainly conjectures based on computer calculations. Some results (see Theorem 4.1 below) were established in previous works, but the aim of this exposition is to point out several mathematical statements remaining largely mysterious and longing for clearness.
 \vskip 0.5 cm
\noindent{\bf{2. Continued fractions}} 
 \par Concerning continued fractions in the area of function fields, we use classical notation, as it can be found for instance in \cite[p.~3-8]{L1}. Throughout the paper we are dealing with finite sequences (or words), consequently we recall the following notation on sequences in $\F_q[T]$. Let $W=w_1,w_2,\ldots,w_n$ be such a finite sequence, then we set $\vert W\vert =n$ for the length of the word $W$. If we have two words $W_1$ and $W_2$, then $W_1,W_2$ denotes the word obtained by concatenation. 
\newline As usual, we denote by $[W]=[w_1,\dots,w_n]\in \F_q(T)$ the finite continued fraction $w_1+1/(w_2+1/(\dots ))$. In this formula the $w_i$, called the partial quotients, are non-constant polynomials. Still, we will also use the same notation if the $w_i$ are constant and the resulting quantity is in $\F_q$. However in this last case, by writing $[w_1,w_2,\dots,w_n]$ we assume that this quantity is well defined in $\F_q$, i.e. $w_n\neq 0,[w_{n-1},w_n]\neq 0,\dots,[w_2,\dots,w_n]\neq 0$.
\newline For $n\geq 0$, a continuant $X_n$ is a polynomial, in the $n$ variables $x_1,x_2,\dots,x_n$, defined recursively by $X_0=1$, $X_1=x_1$ and $X_k=x_kX_{k-1}+X_{k-2}$ for $2\leq k\leq n$. We use the notation $\langle W\rangle$ for the continuant built from $W=w_1,w_2,\dots,w_n$. In the sequel the $w_i's$ are in $\F_q[T]$, then the degree in $T$ of $\langle W\rangle$ is clearly equal to the sum of the degrees in $T$ of the $w_i's$. 
\newline We denote by $W'$ (resp. $W''$) the word obtained from $W$ by removing the first (resp. last) letter of $W$. Hence, we recall that we have $[W]=\langle W\rangle/\langle W'\rangle$. We let $W^*=w_n,w_{n-1},\ldots,w_1$, be the word $W$ written in reverse order. We have $\langle W^*\rangle$=$\langle W\rangle$ and also $[W^*]=\langle W\rangle/\langle W''\rangle$. 
\newline Moreover, if $y\in \F_q^*$, then we define $y\cdot W$ as the following sequence$$y\cdot W = y w_1, y^{-1}w_2,\ldots, y^{(-1)^{n-1}}w_n.$$
Then, it is also known that $\langle y\cdot W\rangle=y\langle W\rangle$ (resp. $=\langle W\rangle$) if  $|W|$ is odd (resp. if $|W|$ is even) and $[y\cdot W]=y[W]$.
\par If $\alpha \in \F(q)$ is an infinite continued fraction, $\alpha=[a_1,a_2,\dots,a_n,\dots]$, we set $x_n=\langle a_1,a_2,\dots,a_n\rangle$ and $y_n=\langle a_2,\dots,a_n\rangle$. In this way, we have $x_n/y_n=[a_1,a_2,\dots,a_n]$, with $x_1=a_1$, $y_1=1$ and by convention $x_0=1$, $y_0=0$. We introduce $\alpha_{n+1}=[a_{n+1},a_{n+2},\dots]$ as the tail of the expansion or the complete quotient ($\alpha_1=\alpha$), and we have 
$$\alpha=(x_{n}\alpha_{n+1}+x_{n-1})/(y_{n}\alpha_{n+1}+y_{n-1})\quad \text{ for }\quad n\geq 1.$$
We recall the following general result (see \cite[p.~332]{L2}).
\newline {\bf{Proposition  2.1. }}{\emph{Let $p$ be a prime number, $q=p^s$ and $r=p^t$ with $s,t\geq 1$. Let $\ell\geq 1$ be an integer and $(a_1,a_2,\dots,a_\ell)\in (\F_q[T])^l$, with $\deg(a_i)>0$ for $1\leq i\leq \ell$. Let $(P,Q)\in (\F_q[T])^2$ with $\deg(Q)<\deg(P)<r$. Then there exists a unique infinite continued fraction $\alpha$ in $\F(q)^{+}$  satisfying
$$(*)\quad \alpha=[a_1,a_2,\dots,a_\ell,\alpha_{\ell+1}]\quad \text{and }\quad \alpha^r=P\alpha_{\ell+1}+Q.$$ This element $\alpha$ is the unique root in  $\F(q)^{+}$ of the following algebraic equation:
$$(**)\quad y_\ell X^{r+1}-x_\ell X^r+(Py_{\ell-1}-Qy_\ell)X-Px_{\ell-1}+Qx_\ell=0.$$}}

\noindent{\bf{3. $P_k$-continued fractions}} 
\par In this note, to avoid unnecessary sophistication in a first stage, the base field will simply be the finite field $\F_p$ having $p$ elements. Moreover, for the subject treated here we need have odd characteristic. Hence, here $p$ is a prime number with $p>2$. Actually, in many areas concerning power series over a finite field, both cases, even and odd characteristic, must be considered separately.
\par Throughout this note $p$ being an odd prime number, $k$ is an integer with $1\leq k<p/2$. We shall now describe several numbers and polynomials which appeared in earlier works. These were introduced in \cite{L2}, to where the reader is invited to refer for the proof of the properties stated below. 
\newline We set $P_k(T)=(T^2-1)^k$. Then we have the following euclidean division :
$$T^p=A_1P_k+R_k.$$
The polynomial $A_1$ is the integer part of the rational $T^p/P_k$, denoted $A_1=[T^p/P_k]$ and $R_k$ is the remainder. Note that $\deg(A_1)=p-2k$ and $A_1=T$ if $k=(p-1)/2$. We define $$\omega_k=(-1)^{k-1}\prod_{1\leq i\leq k}(1-1/2i)\in \F_p^*.$$
In previous works another polynomial closely related to $R_k$ was used. Defining $Q_k=(2k\omega_k)^{-1}R_k$, it was proved that $Q_k$ can also be defined by 
$$Q_k(T)=\int_{0}^{T}(x^2-1)^{k-1}dx=\sum_{0\leq i\leq k-1}(-1)^{k-1-i}\binom{k-1}{i}(2i+1)^{-1}T^{2i+1}.$$ 
These formulas are the key to obtain the following : we have the remarkable continued fraction expansion $$T^p/P_k=[A_1,w_1T,w_2T,\dots,w_{2k}T],\eqno{(1)}$$ 
where the $w_i\in \F_p^*$ are defined by : $w_1=(2k-1)(2k\omega_k)^{-1}$ and
$$w_{i+1}w_i=(2k-2i-1)(2k-2i+1)(i(2k-i))^{-1} \qquad \text{ for
}\quad 1\leq i\leq 2k-1.$$
Note that the $w_i$ are rational numbers which, for all primes $p$ with the condition $2k<p$, exist by reduction in $\F_p^*$. We also have
$$w_{2k+1-i}=-w_i  \quad \text{ for }\quad 1\leq i\leq 2k.\eqno{(2)}$$
In the sequel, we consider the finite word $W=w_1T,w_2T,\dots,w_{2k}T$. Hence, with our notation, $(1)$ implies $T^p/P_k=A_1+1/[W]$, and we can write $$P_k/R_k=[w_1T,w_2T,\dots,w_{2k}T]=[W].\eqno{(3)}$$ 
\par Let us do a thorough investigation of the extremal case: $k=(p-1)/2$.
Then, we clearly have $A_1=[T^p/P_k]=T$ and consequently
$$R_k=T^p-T(T^2-1)^k=T^p-T(T^{2k}-kT^{2k-2}+\cdots)=kT^{2k-1}+\cdots$$
Hence, we get $[P_k/R_k]=(1/k)T=w_1T$. Since $2k=-1$, we have $w_1=1/k=-2$. We obtain $\omega_k=(2k-1)/(2kw_1)=-2/(-w_1)=-1$. But we also have  $w_{i+1}w_i=-4$ and it follows that $W=-2T,2T,-2T,\cdots,2T$. See
\cite[Corollary 4.2]{L2} with a different notation. 
\par Applying Proposition 2.1, we consider particular algebraic continued fractions in $\F(p)$ defined as follows. 
\newline {\bf{Definition 3.1.}}{\emph{ An infinite continued fraction $\alpha\in \F(p)$ is a $P_k$-expansion if there exist an integer $\ell \geq 1$, a $\ell$-tuple, $(\lambda_1,\lambda_2,\dots,\lambda_\ell)\in (\F_p^*)^{\ell}$ and $(\epsilon_1,\epsilon_2)\in \F_p^*\times \F_p$ such that 
$$\alpha=[\lambda_1T,\lambda_2T,\dots,\lambda_{\ell}T,\alpha_{\ell+1}]\quad \text{and }\quad \alpha^p=\epsilon_1P_k\alpha_{\ell+1}+\epsilon_2R_k.$$}}
\newline For given integers $k$, $p$ and $\ell$, chosen as above,  a $P_k$-expansion is defined by a $(\ell+2)$-tuple: $\Lambda=(\lambda_1,\lambda_2,\dots,\lambda_{\ell},\epsilon_1,\epsilon_2)$. Moreover, this continued fraction satisfies the algebraic equation $(**)$, introduced in Proposition 2.1, with $(P,Q)=(\epsilon_1P_k,\epsilon_2R_k)$, whose four coefficients in $\F_p[T]$ are only depending on $k,p,\ell$ and  $\Lambda$.
\par These particular continued fraction expansions, introduced in \cite{L2}, have been considered in several articles and they were particularly studied in a more general setting in \cite{L4}. In previous publications, we used in this definition the polynomial $Q_k$ instead of $R_k$. This has no consequence, since both polynomials are proportional. Originally, a particular and simple continued fraction expansion in $\F(13)$, solution of an algebraic equation of degree four, was observed by Mills and Robbins \cite{MR}. In order to understand this particular and remarkable pattern, we developed this notion of $P_k$-expansion and so we could show that Mills and Robbins example was a particular  case, defined as above with $p=13$, $k=4$ and $\ell=6$. It implied its full description (see \cite{L3} and also \cite{L5,AL} for a generalization concerning Robbins' quartic). The possibility of describing the sequence of partial quotients for an arbitrary $P_k$-expansion is yet out of reach. However, using a technical lemma stated below, we will see in the next section that the sequence of partial quotients for certain $P_k$-expansions can be explicitly given. This lemma is the following (see the origin in \cite[p.~336 and p.343]{L2}). 
\newline {\bf{Lemma 3.2. }}{\emph{ Let $A\in \F_p[T]$, $\delta \in \F_p^*$ and $X\in \F(p)$. Then we have $$[A+\delta R_kP_k^{-1},X]=[A, \delta^{-1}\cdot W,XP_k^{-2}+\delta^{-1}R_kP_k^{-1}].$$ }}
\newline Proof: According to $(3)$, we have $P_k/R_k=[w_1T,w_2T,\dots,w_{2k}T]=[W]$. By $(2)$, we also have $W^*=-1\cdot W$. We observe the following links between the polynomials $P_k,R_k$ and the continuants built from $W$. We have
$$P_k/R_k=[W]=\langle W\rangle/\langle W'\rangle.$$ 
We note that $\langle W\rangle=\langle w_1T,w_2T,\dots,w_{2k}T\rangle$ and $P_k=(T^2-1)^k$ have the same degree in $T$ equal to $2k$. Besides, the constant term of $P_k$ is $(-1)^k$ while, since $W$ has even length, the constant term of $\langle W\rangle$ is $1$. Consequently we obtain
$$\langle W\rangle= (-1)^kP_k\qquad \text{ and
} \quad \langle W'\rangle= (-1)^kR_k.\quad \eqno{(4)}$$
 Hence, we can write $$A+\delta R_k/P_k=A+1/(\delta^{-1}[W])=A+1/[\delta^{-1}\cdot W]=[A,\delta^{-1}\cdot W].$$
We set $U=A,\delta^{-1}\cdot W=u_1,u_2,\dots,u_n$. Here we have $n=2k+1$, $u_1=A$ and $u_{i+1}=\delta^{(-1)^i}w_iT$ for $1\leq i\leq 2k$. We set $x_n=\langle U\rangle$ and $y_n=\langle U'\rangle$, so that we have $[U]=x_n/y_n$. Let us consider an arbitrary element $b\in \F(p)$. Since $x_ny_{n-1}-y_nx_{n-1}=(-1)^n$, we can write 
$$[U,b]-[U]=(x_nb+x_{n-1})/(y_nb+y_{n-1})-x_n/y_n=(-1)^{n-1}(y_n(y_nb+y_{n-1}))^{-1}.$$
 Then we have $y_n=\langle U'\rangle=\langle \delta^{-1}\cdot W \rangle=\langle W\rangle$, since $|W|$, the length of $W$, is even. But also $y_{n-1}=\langle (U')''\rangle=\langle \delta^{-1}\cdot W'' \rangle=\delta^{-1}\langle W''\rangle$, since $|W''|$ is odd. Due to the properties of continuants and by (2), we have :
$$\langle W''\rangle=\langle w_1T,w_2T,\dots,w_{2k-1}T\rangle= \langle  -w_{2k}T,\dots,-w_2T\rangle=-\langle W'\rangle.$$
Hence, by (4), we have
$$ y_n=(1)^kP_k \quad \text{ and }\quad y_{n-1}=(-1)^{k+1}\delta^{-1}R_k.$$
Therefore we obtain
$$[U,b]-[U]=(P_k^2b-\delta^{-1}P_kR_k)^{-1}.$$
Choose $b=XP_k^{-2}+\delta^{-1}R_kP_k^{-1}$ or equivalently $X=P_k^2b-\delta^{-1}P_kR_k$. The last formula gives $[U,b]=[U]+X^{-1}$ and therefrom we obtain 
$$[A,\delta^{-1}\cdot W,XP_k^{-2}+\delta^{-1}R_kP_k^{-1}]=[U]+X^{-1}=[A+\delta R_kP_k^{-1},X].$$
So the proof is complete.

\noindent{\bf{4. Perfect $P_k$-continued fractions}}
\par Our goal is to describe explicitly certain $P_k$-expansions. This will be done in two different cases, each corresponding to a particular choice of  $\Lambda$. Because of some change of notation and to highlight the similarity between both cases, we repeat the results concerning the first case, which have already appeared in previous publications. To describe the partial quotients in these continued fractions, we need to introduce in $\F_p[T]$ two sequences of polynomials  $(A_{n})_{n\geq 0}$ and $(B_{n})_{n\geq 0}$ as follows. The first one is defined by
$$A_{0}=T \quad \text{ and recursively }\quad  A_{n+1}=[A_{n}^p/P_k] \quad \text{ for }\quad n\geq 0.$$ Here the brackets denote the integral (i.e. polynomial) part of the rational function. Note that, in agreement to the beginning of the previous section, we have $A_1=[T^p/P_k]$. While the second sequence is defined by
$$B_{0}=A_{0}=T \quad \text{ and }\quad B_{1}=A_{1}=[T^p/P_k]$$ 
and recursively $$B_{n+1}=B_{n}^pP_k^{(-1)^{n+1}}\quad \text{ for }\quad n\geq 1.$$ 
We are particularly interested in the degrees of these polynomials. We set $u_n=\deg(A_{n})$ and $v_n=\deg(B_{n})$. From the recursive definition of these polynomials, 
we get $u_0=v_0=1$ and also 
$$u_{n+1}=pu_n-2k \quad \text{ and }\quad v_{n+1}=pv_n+2k(-1)^{n+1}\quad \text{ for }\quad n\geq 0.$$ 
Note that the sequence $(u_n)_{n\geq 0}$ is constant if $2k=p-1$, then we have $A_n=A_0=T$ for $n\geq 0$. Otherwise, 
both sequences  $(u_n)_{n\geq 0}$ and $(v_n)_{n\geq 0}$ are strictly increasing. 
\newline The first case is described in the following theorem.
\newline {\bf{Theorem 4.1. }}{\emph{ Let $p$ be an odd prime and $k,\ell$ integers chosen as above. Let $\alpha\in \F(p)$ be a $P_k$-expansion, depending as above on a $(\ell+2)$-tuple $\Lambda$. We assume that $\epsilon_2\neq 0$ and $\Lambda$ satisfies the following equality:
 $$[\lambda_\ell,\lambda_{\ell-1},\dots,\lambda_2,\lambda_1-\epsilon_2]=\epsilon_1/(\epsilon_2\omega_k).\eqno{\mathcal{H}(A)}$$
Then there exist a sequence $(\lambda_n)_{n\geq 1}$ in $\F_p^*$ and a sequence $(i(n))_{n\geq 1}$ in $\N$, such that
$$a_n=\lambda_n A_{i(n)}\quad \text{ for }\quad n\geq 1.$$}}

\par The proof of this theorem appears in \cite[p.~1111]{L3}, but also in  \cite[p.~256-257]{L4}, in a much larger context. See the remark in \cite[p.~257]{L4}, concerning a mistake in the writing of a formula in \cite{L3}. This proof is derived from Lemma 3.2 presented in Section 3. Both sequences $(\lambda_n)_{n\geq 1}$ and $(i(n))_{n\geq 1}$ have been described in \cite{L3} and \cite{L4}, sometimes with different notation. We recall here below the description of these sequences, when the base field is prime, which is the only case we consider here. 
\newline For $n\geq 1$, we set $f(n)=(2k+1)n+\ell -2k$. At first, according to $\mathcal{H}(A)$, we can define in $\F_p^*$ : $z_1=\lambda_1-\epsilon_2$ and $z_i=\lambda_i+z_{i-1}^{-1}$ for $2\leq i\leq \ell$. Then we consider the sequence $(z_{n})_{n\geq 1}$ in $\mathbb{F}_{p}^{\ast }$ defined, from the initial 
values $z_1,z_2,\dots,z_{\ell}$, by the recursive formulas
$$z_{f(n)+i}=\theta _{i}\epsilon _{1}^{(-1)^{n+i}}z_n^{(-1)^{i}}\quad \text{for}\ n\geq 1\quad \text{and}\quad 0\leq i\leq 2k,\eqno{(A_1)}$$ 
where 
$$\theta _{0}=-\omega_k\quad \text{and}\quad \theta_i=-iw_i/(2k-2i+1)\quad \text{for}\quad 1\leq i\leq 2k.$$ 
 Then the sequence $(\lambda
_{n})_{n\geq 1}$ in $\mathbb{F}_{p}^{\ast }$ is defined recursively,
from the initial values $\lambda _{1},\lambda _{2},\ldots ,\lambda _{\ell }$, by the formulas
$$\lambda _{f(n)}=\epsilon _{1}^{(-1)^{n}}\lambda
_{n} \quad \text{and}\quad \lambda _{f(n)+i}=\epsilon _{1}^{(-1)^{n+i}}w_iz_n^{(-1)^i},\eqno{(A_2)}$$
for $n\geq 1$ and for $1\leq i\leq 2k$.
\newline It is worth to mention that the complexity of this sequence $(\lambda
_{n})_{n\geq 1}$ in $\mathbb{F}_{p}$ has been studied. Indeed, it was proved to be $(2k+1)$-automatic. For a full account on this matter, in the most general setting, the reader is advised to refer to \cite[Section 5]{LY2}.

\par Concerning the sequence $(i(n))_{n\geq 1}$, we have the following description: 
$$i(n)=0 \quad \text{ if} \quad n\notin f(\N^*)\quad \text{and}\quad i(f(n))=i(n)+1 \quad \text{for}\quad n\geq 1.$$
\par We want to compute the irrationality measure for a continued fraction described in Theorem 4.1. We only need to know the sequence $(i(n))_{n\geq 1}$ and the sequence $(u_n=\deg(A_n))_{n\geq 0}$. The proof of the following corollary is based on another equivalent description of this sequence $(i(n))_{n\geq 1}$ (See \cite[p.~143]{AL}). In the sequel, we will use the following notation. For $n\geq 0$, if we have the word $w,w,\dots,w$ of length $n$, then we denote it shortly by $w^{[n]}$ with $w^{[0]}=\emptyset$. In the same way $W^{[n]}$ denotes the word $W,W,\dots,W$ where $W$ is a finite word repeated $n$ times and $W^{[0]}=\emptyset$. Let  $(I_n)_{n\geq 0}$ be the sequence of finite words of integers defined recursively by
$$I_0=0\quad \text{ and }\quad I_n=n,I_0^{[2k]},I_1^{[2k]},\cdots,I_{n-1}^{[2k]}\quad \text{ for }n\geq 1.$$
Then the sequence $I=(i(n))_{n\geq 1}$ in $\N$, introduced in Theorem 4.1, is given by the infinite word:
$$I=I_0^{[\ell]},I_1^{[\ell]},I_2^{[\ell]},\cdots,I_n^{[\ell]},\cdots $$
By mean of the last formula, on the irrationality measure, given in the introduction, we get the following result (see \cite[p.~148]{AL}).
\newline {\bf{Corollary 4.2. }}{\emph{ Let $\alpha$ be a $P_k$-expansion, depending  on the $(\ell+2)$-tuple $\Lambda$ satisfying $\mathcal{H}(A)$, then we have
$$\nu(\alpha)=2+(p-2k-1)/\ell.$$}}
\par Before going on to a second kind of $P_k$-expansion, we need to point at an exceptional element in $\F(p)$ for all $p\geq 3$ which is quadratic but also a perfect $P_k$-expansion of the type described in Theorem 4.1. This element appears in the extremal case $k=(p-1)/2$. We have observed that, in this case, the sequence $(A_n)_{n\geq 0}$ is constant and $A_n=T$ for $n\geq 0$. Consequently, the elements described in theorem 4.1 will have all partial quotients of the form $a_n=\lambda_nT$ (note that, in agreement to Corollary 4.2, the irrationality measure of the continued fraction will be equal to 2). In these cases the sequence $(\lambda_n)_{n\geq 1}$ will not generally be periodic and this is why these non-quadratic examples were brought to light in \cite{MR}. However at the bottom of these examples lies a universal element which is the formal golden mean. In $\F(p)$, for $p\geq 3$, let us consider the following infinite continued fractions
$$\phi(T)=[T,T,\cdots,T,\cdots] \quad \text{and}\quad \rho(T)=[-2T,2T,\cdots,-2T,2T,\cdots].$$ 
Note that we have $\rho(T)=(u/2)\phi(uT)$ with $u^2=-4$. Elementary computations (see \cite[p.~331-332]{L2} and \cite[p.~267-268]{LY1})  show that we have 
$$\rho^p=(-1)^{\ell}P_k\rho_{\ell+1}-R_k \quad \text{for}\quad k=(p-1)/2 \quad \text{and}\quad \ell \geq 1.$$
Since $\omega_k=-1$, we observe that $\mathcal{H}(A)$ reduces to $$[(-1)^{\ell}2,(-1)^{\ell-1}2,\cdots,-2+1]=(-1)^{\ell},$$
which is easily verified by induction. This element $\rho$ will appear again in the last section (Proposition 5.1).
\par If we have reported here Theorem 4.1, this is due to its somehow mysterious proof and also to the apparent closeness with the following conjecture. There is a second case, corresponding to a different choice of $\Lambda$,  where a $P_k$-expansion can be described but only partially and conjecturally.
\newline {\bf{Conjecture 4.3. }}{\emph{ Let $p$ be an odd prime and $k,\ell$ integers chosen as above. Let $\alpha\in \F(p)$ be a $P_k$-expansion, depending as above on a $(\ell+2)$-tuple $\Lambda$. We assume that $\Lambda$ satisfies the following equality:
 $$[\lambda_\ell,\lambda_{\ell-1},\dots,\lambda_1-\epsilon_2]=0.\eqno{\mathcal{H}(B)}$$
Then there exist a sequence $(\lambda_n)_{n\geq 1}$ in $\F_p^*$ and a sequence $(j(n))_{n\geq 1}$ in $\N$, such that
$$a_n=\lambda_n B_{j(n)}\quad \text{ for }\quad n\geq 1.$$}}
\par It may be worth to notice that condition $\mathcal{H}(B)$ (and also $\mathcal{H}(A)$) can be stated in a different way. Indeed $\mathcal{H}(B)$ is simply equivalent to : $\epsilon_2=[\lambda_1,\cdots,\lambda_{\ell-1},\lambda_{\ell}]$. This conjecture results from broad computer observations, letting the parameters $p,k$ and $\ell$ vary. We also have a conjectural description of the sequence $(j(n))_{n\geq 1}$ as follows.
\newline Let  $(J_n)_{n\geq 1}$ be the sequence of finite words of integers defined recursively by
$$J_0=0\quad \text{ and }\quad J_{n}=n,J_{n-1}^{[2k-1]},n-1\quad \text{ for }n\geq 1.$$
Then the sequence $J=(j(n))_{n\geq 1}$ in $\N$ is given by the infinite word:
$$J=J_0^{[\ell-1]},0,J_1^{[\ell-1]},1,\cdots,J_n^{[\ell-1]},n,\cdots $$
 The irrationality measure, for a continued fraction described in Conjecture 4.3, only depends on the sequence $(j(n))_{n\geq 1}$ and the sequence $(v_n=\deg(B_n))_{n\geq 0}$. Accordingly, using arguments similar to the ones given in \cite[p.~150-151]{AL}), we obtain the following.
\newline {\bf{Conjecture 4.4. }}{\emph{ Let $\alpha$ be a $P_k$-expansion, depending  on the $(\ell+2)$-tuple $\Lambda$ satisfying $\mathcal{H}(B)$, then we have
$$\nu(\alpha)=2+(p-2k+1)(p-1)/(\ell(p+1)-2k).$$}}
The description of the sequence $(\lambda_n)_{n\geq 1}$ in $\F_p^*$ is in general yet out of reach. Note that in the simplest case $\ell=1$, it can be proved that $J=\N$ and $a_n=\lambda_nB_{n-1}$ where $\lambda_{2n}=\lambda_1\epsilon_1^{-1}$ and $\lambda_{2n+1}=\lambda_1$ for $n\geq 1$. Moreover, in this case, the irrationality measure is equal to $p+1$ which is the maximal possible value. In \cite[p.~20-21]{L1}, the sequence of partial quotients has also been fully described in the simple case: $k=1$ and $\ell=2$.
\par The reader may wonder why is a $P_k$-expansion particular in the two cases presented above. In these cases we say that the expansion is perfect of type A or of type B. First, a computer observation, outside these cases, shows a certain irregularity in the sequence of partial quotients. Moreover the specificity of these perfect expansions of type A and B is pointed out in the next and last section. Finally, the possible connection between these two cases is an open question (regarding the similarity between $\mathcal{H}(A)$ and  $\mathcal{H}(B)$).
\par Before concluding this section, we must make a remark on the origin of these particular continued fractions. As pointed out above, the $P_k$-expansions appeared from the study of Mills and Robbins quartic over $\F_{13}$ \cite[p.~403-404]{MR}. This quartic equation has been generalized for all primes $p\geq 5$ in \cite{L5}. The solution has a continued fraction expansion with a different pattern according to the residue modulo 3 of the prime $p$. In the case $p\equiv 1 \mod 3$ (particularly for $p=13$), the pattern of this continued fraction is perfect of type A. In the case $p\equiv 2 \mod 3$, a different pattern, described in \cite{AL}, appears and the partial quotients are all proportional to elements of the sequence $(B_n)_{n\geq 0}$. However the continued fraction is not of the same type as in Conjecture 4.3 : it is a hyperquadratic element of order 2. This has given rise to the following generalization concerning this sequence $(B_n)_{n\geq 0}$, in connection with hyperquadratic elements of higher order.
\par Returning to Proposition 2.1, let us take $r=p^n$ with $n\geq 1$. Let us make the following choice for the pair $(P,Q)\in (\F_p[T])^2$: 
$$P=\epsilon_1P_k^{(p^n+(-1)^{n-1})/(p+1)}\quad \text{ and }\quad Q=\epsilon_2R_k^{p^{n-1}}\quad \text{where }\quad \epsilon_1\neq 0,\epsilon_2\in \F_p.$$ 
We consider the continued fraction $\alpha$ in $\F(p)$ so defined by $(*)$, which is a hyperquadratic element of order $n$. We observed the following. If the $(\ell+2)$-tuple $\Lambda=(a_1,\cdots,a_{\ell},\epsilon_1,\epsilon_2)$ is well chosen, then all partial quotients for $\alpha$ appear to be proportional to certain polynomials $B_n$ as in Conjecture 4.3 (corresponding to the case $n=1$). Again this observation is just based on broad computer calculations. Amazingly, this is what happens for the solution of the generalized quartic when $p\equiv 2 \mod 3$, in this case we have $n=2$ and we will say that the expansion is perfect of type B and order $2$. In \cite[p.~141-143]{AL}), we have described a choice of $\Lambda$ to obtain the above conjecture in the case $n=2$ and the description of the corresponding analogue of the sequence $(j(n))_{n\geq 1}$. 
\par Here below, to support the observation which has just been made, we have collected a few examples of presumably perfect $P_k$-expansions of type B and of order $3$. To allow easier computer calculations, we have only considered the values $p=3$ and $k=1$.
\newline {\bf{ Conjectural examples 4.5. }}{\emph{ ($n=3$, $p=3$, $k=1$ and $\epsilon_1=\epsilon_2=1$).
\newline Let $\ell \geq 1$ be an integer. Let $\alpha=[a_1,a_2,\cdots,a_n,\cdots]\in \F(3)$ be the expansion defined by
$$ \alpha=[a_1,a_2,\dots,a_l,\alpha_{\ell+1}]\quad \text{and }\quad \alpha^{27}=(T^2-1)^7\alpha_{\ell+1}+T^9$$
where we have $\ell=3$ and $(a_1,a_2,a_3)=(2T,T,T)$
\newline or $\ell=4$ and $(a_1,a_2,a_3,a_4)=(T,2T,2T,T)$
\newline or $\ell=4$ and $(a_1,a_2,a_3,a_4)=(2T,T^5+2T^3,T,T)$
\newline or $\ell=5$ and $(a_1,a_2,a_3,a_4,a_5)=(2T,T,T,T,T^5+2T^3)$. 
\newline Then, in all these cases, we have observed that there exist $\lambda_n\in \F_3^*$ and $j(n)\in\N$ such that we have $a_n=\lambda_nB_{j(n)}$ up to a certain rank. We conjecture that the same holds for all $n\geq 1$.}}
\newline Let us make a comment about this last statement. There are several ways to obtain the first partial quotients of an algebraic element $\alpha$ in $\F(q)$ (see \cite[p.~34]{L5}). A natural way is to obtain the beginning of the continued fraction expansion from a rational approximation to $\alpha$. In the present case (this is valid for all hyperquadratic elements), we can write the algebraic equation $(**)$ satisfied by $\alpha$ in the following way : $\alpha=(a\alpha^r+b)/(c\alpha^r+d)=g(\alpha)$. Then we build rational approximations to $\alpha$ , starting from $R_1=x_\ell/y_\ell$, by the recursive formula $R_{n+1}=g(R_n)$ for $n\geq 1$. In the cases indicated above, these rational numbers have a special form, showing a particular expansion with partial quotients as wanted. We guess a proof of our conjecture could be derived from the study of the numerators and denominators of $R_n$.
\par At last, we remark the existence of another different generalization for $P_k$-expansions, also leading to hyperquadratic elements of higher order, which was presented in \cite{L4}. In the next and last section, we return to $P_k$-expansions of order $1$ (i.e. $r=p$).
\vskip 0.5 cm

\noindent{\bf{5. Hyperquadratic elements having a low algebraic degree}} 
\par In this section, we introduce a family of polynomials in the variable $X$ with coefficients in $\F_p[T]$. These polynomials $H$ have the particular form of hyperquadratic type :
$$H(X)=uX^{p+1}+vX^p+wX+z \quad \text{ with }\quad u,v,w,z \in \F_p[T].$$ 
The four coefficients are depending on parameters coming from the finite word $W$ introduced in Section 3. Moreover $H(X)=0$ will have solutions in $\F(p)$ having a continued fraction expansion of the type discussed above. Under a particular choice of these parameters, $H$ will be reducible and consequently these solutions will have a particular algebraic degree smaller than $p+1$. The consideration of these polynomials comes from a quartic equation, generalizing Mills and Robbins example, introduced in \cite{L5}. To describe these coefficients, we use some more notation appearing there \cite[p.~30-31]{L5}.
\newline As above $p$ is an odd prime number and $k$ an integer with $1\leq k<p/2$ and $W=w_1T,w_2T,\cdots,w_{2k}T$ is the finite word introduced in Section 3.
\newline Let $n,t$ and $m$ be integers with $1\leq n\leq t\leq m\leq 2k$. Then we introduce the continuant :
$$K_{n,m}=\langle w_nT,w_{n+1}T,\cdots,w_mT\rangle \in \F_p[T].$$ 
By convention we extend this notation with $K_{n,n-1}=1$ and $K_{n,n-2}=0$. From general properties of continuants (see \cite[p.~7]{L1}), we have the following formula 
$$K_{n,m}K_{t,m-1}-K_{n,m-1}K_{t,m}=(-1)^{m-t}K_{n,t-2}.\eqno{(5)}$$
Since, by $(2)$, we have $W^*=-1\cdot W$, $\langle A^*\rangle=\langle A\rangle$ and $\langle -1\cdot A\rangle=(-1)^{\vert A\vert}\langle A\rangle$ for any finite word $A$, we also get
$$K_{n,m}=(-1)^{m-n+1}K_{2k+1-m, 2k+1-n}.\eqno{(6)}$$
Let us consider the vector $V=[p,k,j,\epsilon]$ where $p,k$ are as above, $j$ is an integer with $1\leq j\leq 2k-1$ and $\epsilon \in \F_p^*$. To each such vector $V$, we associate a polynomial $H(V)$ in $\F_p[T][X]$ defined by 
$$H(V)(X)=K_{j+2,2k}X^{p+1}-K_{j+1,2k}X^p+\epsilon (K_{1,j}X+K_{1,j-1}).$$
We have the following proposition.
\newline {\bf{Proposition 5.1. }}{\emph{Let $V=[p,k,j,\epsilon]$ and $H(V) \in\F_p[T][X]$ be defined as above. 
For each choice of $V$, there are two infinite continued fractions $\alpha$ and $\beta$ in $\F(p)$, both $P_k$-expansions, such that we have $H(V)(\alpha)=0$ and  $H(V)(1/\beta)=0$. The first one $\alpha$ is defined by
$$\alpha=[w_{j+1}T,w_{j+2}T,\cdots,w_{2k}T,\alpha_{2k-j+1}]$$
and $$ \alpha^p=\epsilon(-1)^{k+j}(P_k\alpha_{2k-j+1}-R_k).$$
The second one $\beta$ is defined by 
$$\beta=[w_{2k-j+1}T,w_{2k-j+2}T,\cdots,w_{2k}T,\beta_{j+1}]$$
and $$ \beta^p=\epsilon^{-1}(-1)^{k+j}(P_k\beta_{j+1}-R_k).$$
\newline Moreover, there are two different values $\epsilon_A$ and $\epsilon_B$ in $\F_p^*$, corresponding to each triple $(p,k,j)$, such that :
\newline If $\epsilon=\epsilon_A$ then both $P_k$-expansions $\alpha$ and $\beta$ are perfect of type A \newline (i.e. condition $\mathcal{H}(A)$ in Theorem 4.1 is satisfied).
\newline If $\epsilon=\epsilon_B$ then both $P_k$-expansions $\alpha$ and $\beta$ are perfect of type B \newline (i.e. condition $\mathcal{H}(B)$ in Conjecture 4.3 is satisfied).
\newline We have
$$\epsilon_A=(-1)^{k+j+1}[w_{j+1},\cdots,w_{2k},\omega_k]\quad \text{and}\quad \epsilon_B=(-1)^{k+j+1}[w_{j+1},\cdots,w_{2k}].$$
Finally, let $k=(p-1)/2$ and $V=[p,k,j,\epsilon_A]$. We set  $P(X)=X^2+2(-1)^jX+1$, then $P$ divides $H(V)$. We have $\alpha=(-1)^j\rho$ and $P(\alpha)=0$. }}
\newline Proof: According to $(6)$, we can write
$$K_{1,j}=(-1)^{j}K_{2k+1-j, 2k}\quad \text{and}\quad K_{1,j-1}=(-1)^{j-1}K_{2k+2-j, 2k},$$
$$K_{j+1,2k}=(-1)^{j}K_{1,2k-j}\quad \text{and}\quad K_{j+2,2k}=(-1)^{j-1}K_{1,2k-j-1}.$$
Consequently, we get $(-1)^{j-1}H(V)(X)=$
$$K_{1,2k-j-1}X^{p+1}+K_{1,2k-j}X^p+\epsilon (-K_{2k+1-j, 2k}X+K_{2k+2-j, 2k}).$$
and therefore $\epsilon^{-1}(-1)^{j-1}X^{p+1}H(V)(1/X)=$
$$K_{2k+2-j, 2k}X^{p+1}-K_{2k+1-j, 2k}X^p+\epsilon^{-1}(K_{1,2k-j}X+K_{1,2k-j-1}).$$
We set $V^*=[p,k,2k-j,\epsilon^{-1}]$. Hence, we have obtained
$$X^{p+1}H(V)(1/X)=\epsilon (-1)^{j-1}H(V^*)(X).\eqno{(7)}$$
By Proposition 2.1, we know that $\alpha \in \F(p)$ defined by
$$\alpha=[a_1,a_2,\cdots,a_\ell,\alpha_{\ell+1}] \quad \text{and }\quad \alpha^p=P\alpha_{\ell+1}+Q$$
satisfies $I(\alpha)=0$ where we have $I\in \F_p[T][X]$ and
$$I(X)=y_\ell X^{p+1}-x_\ell X^p+(Py_{\ell-1}-Qy_\ell)X-Px_{\ell-1}+Qx_\ell.$$
In order to prove that $H(V)(\alpha)=0$, for the element $\alpha\in \F(p)$ defined in this proposition, we will show that $I=H(V)$ if 
$$(a_1,\cdots,a_\ell,P,Q)=(w_{j+1}T,\cdots,w_{2k}T,\epsilon_1P_k,-\epsilon_1R_k)\quad \text{with}\quad \epsilon_1=(-1)^{k+j}\epsilon.$$
First, recall that, according to $(4)$ in Section 3, we have
$$P_k=(-1)^kK_{1,2k} \quad \text{and}\quad R_k=(-1)^kK_{2,2k}.$$
Consequently, by $(6)$ and $\epsilon_1=(-1)^{k+j}\epsilon$, we get
$$P=\epsilon (-1)^jK_{1,2k} \quad \text{and}\quad Q=\epsilon (-1)^jK_{1,2k-1}.\eqno{(8)}$$
Since $(a_1,\cdots,a_\ell)=(w_{j+1}T,\cdots,w_{2k}T)$, we also have
$$x_\ell=K_{j+1,2k},\quad y_\ell=K_{j+2,2k},\quad x_{\ell-1}=K_{j+1,2k-1},\quad y_{\ell-1}=K_{j+2,2k-1}.\eqno{(9)}$$
Hence we get, from $(8)$, $(9)$ and $(5)$,
$$Py_{\ell-1}-Qy_\ell=\epsilon(-1)^j(K_{1,2k}K_{j+2,2k-1}-K_{1,2k-1}K_{j+2,2k})=\epsilon K_{1,j}$$
and also 
$$Qx_\ell-Px_{\ell-1}=\epsilon(-1)^j(K_{1,2k-1}K_{j+1,2k}-K_{1,2k}K_{j+1,2k-1})=\epsilon K_{1,j-1}.$$
Finally, we obtain the desired outcome
$$I(X)=K_{j+1,2k}X^{p+1}-K_{j+2,2k}X^p+\epsilon K_{1,j}+\epsilon K_{1,j-1}=H(V)(X).$$
Changing $j$ into $2k-j$ and $\epsilon$ into $\epsilon^{-1}$, in the definition for $\alpha$, we get the definition for $\beta$. Since this means changing $V$ into $V^*$, we have just proved that $H(V^*)(\beta)=0$ and according to $(7)$ this gives $H(V)(1/\beta)=0$. 
\newline Let us now consider the $(2k-j+2)$-tuple $\Lambda_{\alpha}$ defining the $P_k$-expansion $\alpha$. We have $\Lambda_{\alpha}=(w_{j+1},\cdots,w_{2k},-\epsilon_2,\epsilon_2)$ where $\epsilon_2=(-1)^{k+j+1}\epsilon$. Hence $\mathcal{H}(A)$ is satisfied if and only if we have $[w_{2k},w_{2k-1}\cdots,w_{j+1}-\epsilon_2]=-1/\omega_k$. This is equivalent to $\epsilon_2=[w_{j+1},\cdots,w_{2k}+1/\omega_k]$. Consequently $\mathcal{H}(A)$ is satisfied by $\Lambda_{\alpha}$ if and only if we have
$$\epsilon=\epsilon_A=(-1)^{j+k+1}[w_{j+1},\cdots,w_{2k},\omega_k].$$
Considering the same $(2k-j+2)$-tuple $\Lambda_{\alpha}$, we see in the same way that $\Lambda_{\alpha}$ satisfies $\mathcal{H}(B)$ if and only if we have
$$\epsilon=\epsilon_B=(-1)^{j+k+1}[w_{j+1},\cdots,w_{2k}].$$
Now, we need to consider $\beta$ instead of $\alpha$. Here the $(j+2)$-tuple $\Lambda_{\beta}$ defining $\beta$ is $\Lambda_{\beta}=(w_{2k-j+1},\cdots,w_{2k},-\epsilon_2,\epsilon_2)$ where $\epsilon_2=(-1)^{k+j+1}\epsilon^{-1}$. In the same way, we obtain that $\Lambda_{\beta}$ satisfies $\mathcal{H}(A)$ if
$$\epsilon=\overline{\epsilon_A}=(-1)^{j+k+1}[w_{2k-j+1},\cdots,w_{2k},\omega_k]^{-1}$$
and also that $\Lambda_{\beta}$ satisfies $\mathcal{H}(B)$ if
$$\epsilon=\overline{\epsilon_B}=(-1)^{j+k+1}[w_{2k-j+1},\cdots,w_{2k}]^{-1}.$$
We only need to prove that $\epsilon_A=\overline{\epsilon_A}$ and  $\epsilon_B=\overline{\epsilon_B}$. Let us compare $\epsilon_B$ and $\overline{\epsilon_B}$. According to $(2)$, we can write $[w_{2k-j+1},\cdots,w_{2k}]=[-w_j,\cdots,-w_1]$. Hence $\epsilon_B=\overline{\epsilon_B}$ is equivalent to 
$$[w_{j+1},\cdots,w_{2k}]+[w_j,\cdots,w_1]^{-1}=0.$$
With our notation on continuants, this can be written as 
$$(K_{j+1,2k}/K_{j+2,2k})(1)+(K_{1,j-1}/K_{1,j})(1)=0.$$
This comes from $K_{1,2k}(1)=(-1)^k(1-1)^k=0$ and the following general formula on continuants: $K_{1,2k}=K_{1,j}K_{j+1,2k}+K_{1,j-1}K_{j+2,2k}$ (see \cite[p~7]{L1}). Turning to $\epsilon_A=\overline{\epsilon_A}$, we see, using the same arguments as above that this is equivalent to 
$$[w_{j+1},\cdots,w_{2k},\omega_k]+[w_j,\cdots,w_1,-\omega_k]^{-1}=0.$$
Using the same general formula on continuants, this is again equivalent to
$$\langle -\omega_k,w_1,w_{2},\cdots,w_{2k},\omega_k\rangle=0.$$
The truth of that is derived from the following two equalities
$$\langle w_{2},\cdots,w_{2k}\rangle=(-1)^{k}\quad \text{and}\quad \langle w_{2},\cdots,w_{2k-1}\rangle=2(-1)^{k+1}\omega_k.$$
The proof of the first one is obtained from $R_k(1)=(-1)^kK_{2,2k}(1)=1$, while the proof for the second one is more mysterious and it is left to the reader. 
\newline We finish the proof of this proposition by considering the extremal case $k=(p-1)/2$. In this particular case, we have seen in Section 3, that we have $W=-2T,2T,-2T,\cdots,2T$ and $\omega_k=-1$. Hence, for all $j$, we get $$\epsilon_A=(-1)^{k+j+1}[2(-1)^{j+1},\cdots,2,-1]=(-1)^k.$$ Consequently, setting $\ell=2k-j$, $\alpha$ is defined by
$$\alpha=[2(-1)^{\ell-1}T,\cdots,2T,\alpha_{\ell+1}] \quad \text{and} \quad \alpha^p=(-1)^{j}(P_k\alpha_{\ell+1}-R_k).$$
Let us consider in $\F(p)$ the infinite continued fraction $\gamma=(-1)^j\alpha$. We have $\gamma^p=(-1)^j\alpha^p$ and $\gamma_{\ell+1}=(-1)^j\alpha_{\ell+1}$ for all $\ell \geq 1$.
Since $(-1)^j=(-1)^{\ell}$, we see that $\gamma$ is defined by
$$\gamma=[-2T,\cdots,2(-1)^{\ell}T,\gamma_{\ell+1}] \quad \text{and} \quad \gamma^p=(-1)^{\ell}P_k\gamma_{\ell+1}-R_k.$$
Comparing to the definition of the quadratic continued fraction $\rho$, derived from the formal golden mean, described in Section 4 and because of the uniqueness in this definition, we obtain $\gamma=\rho$ and $\alpha=(-1)^j\rho$. A basic computation shows that the minimal polynomial of $(-1)^j\rho$ is $P(X)=X^2+2(-1)^{j}TX+1$. Consequently $H$ is a multiple of $P$ and this completes the proof of the proposition.

\par Concerning this polynomial $H(V)$, we are interested in its quality of being or not being reducible over $\F_p(T)$. In other words, we are interested in the exact algebraic degree of $\alpha$ and $\beta$ over $\F_p(T)$. We have checked the reducibility of $H(V)$ by computer calculations. In the vector $V$, we first fix $p$ and $k$. Then, according to formula $(7)$, we observe that the reducibility of $H(V)$ need only be studied for $1\leq j\leq k$. 
\par To illustrate our purpose, let us consider the simplest, and somehow trivial, case : $p=3$. Then, we only have to consider $k=j=1$. Here we have $\omega_1=1/2=-1$ and $W=T,-T$.  There are two polynomials $H$ :
$$H_1(X)=H[3,1,1,1](X)=X^4+TX^3+TX+1$$
and   
$$H_2(X)=H[3,1,1,2](X)=X^4+TX^3-TX-1=(X^2-1)(X^2+TX+1)$$
Note that reducibility of the polynomials $H(V)$ over $\F_p(T)$ has been tested by computer, using Maple programming software. The polynomial $H_1$ is irreducible over $\F_3(T)$. According to Proposition 5.1, in both cases, $\alpha$ is defined by $\alpha=[-T,\alpha_2]$ and $\alpha^3=\epsilon((T^2-1)\alpha_2-T)$ and we have $\alpha=\beta$ (hence $\alpha$ and $1/\alpha$ are solutions in $\F(3)$ of $H$). Moreover, we get $\epsilon_A=(-1)^3[2,-1]=-1$ and $\epsilon_B=(-1)^3[2]=1$. If $\epsilon=1$, the solution $\alpha\in \F(3)$ of $H_1(X)=0$ is algebraic of  degree 4, we have $\alpha=-[B_0,B_1,\cdots,B_n,\cdots]$ and $\nu(\alpha)=4$ (see the comment after Conjecture 4.4). If $\epsilon=-1$, the solution $\alpha$ of $H_2(X)=0$ satisfies $\alpha^2+T\alpha+1=0$ and therefrom we get $\alpha=[-T,T,\cdots,-T,T,\cdots]$. Indeed $\alpha=u\phi(uT)$ where $\phi(T)=[T,T,T,\dots,T,\dots]$ is the formal golden mean and $u^2=-1$. 

\par Given the prime $p$, an easy computation show that there are $(p^2-1)/8$ pairs $(j,k)$ with $1\leq j\leq k\leq (p-1)/2$. Hence, for a given odd prime number $p$, there are $(p^2-1)(p-1)/8$ polynomials $H(V)$ to be considered. We have investigated the reducibility of these polynomials in all the possible cases, for the small primes $p$. Based on this observation, we can make the first and main conjecture. Note that the three results presented below are clearly true, just by direct computation, with limitation on the size of the prime $p$.
\newline {\bf{Conjecture 5.2. }}{\emph{Let $p$ be an odd prime number. Let $V=[p,k,j,\epsilon]$ and $H(V)$ be defined as above. If $\epsilon\neq \epsilon_A$ and $\epsilon\neq \epsilon_B$, then the polynomial $H(V)$ is irreducible over $\F_p(T)$.}} 

\par The irreducibility of $H$, outside the two values for $\epsilon$ given in the proposition, may be surprising. Another surprise is that the reducibility may sometimes happen, for a particular value of $\epsilon$, but only if the triple $(p,k,j)$ is well chosen. A last surprise is the following : if the reducibility appears  for a certain triple $(p,k,j)$ and a particular choice of $\epsilon$, then this is most of the times if $\epsilon=\epsilon_A$, and sometimes if $\epsilon=\epsilon_A$ or $\epsilon=\epsilon_B$, but apparently never if $\epsilon=\epsilon_B$ alone. To illustrate this, we give below three tables where we have collected the data obtained by computer for $p=7$, $p=11$ and $p=13$. Each cell $(j,k)$ with $j\leq k$, corresponds to $p-1$ tested polynomials. In each cell, we have written Ir. if the polynomials $H$ are irreducible for all $\epsilon$. If $H$ is reducible (for $\epsilon=\epsilon_A$ or $\epsilon=\epsilon_B$), we indicate the degrees of its factors ( For instance, $A:5^2.2$ means 2 factors of degree 5 and 1 of degree 2, with $\epsilon=\epsilon_A$ and $B:8.4$ means 1 factor of degree 8 and 1 of degree 4, with $\epsilon=\epsilon_B$. Note that the sum of these degrees is $p+1$).
\vskip 0.5 cm
\centerline{Reducibility of $H$: $p=7$}
\begin{center}
   \begin{tabular}{| c | c | c | c | }
     \hline
      & 1 & 2 & 3 \\ \hline
     3 & $A:3^2.2$ & $A:3^2.2$ & $A:2^3.1^2$ \\ \hline
     2 & $A:4^2$ & $A:3^2.1^2$ & \\     
      &  & $B:4^2$ & \\ \hline
         1 & $A:6.1^2$ &  &  \\
                   \hline           
   \end{tabular}
   \end{center}
   \vskip 0.5 cm
\centerline{Reducibility of $H$: $p=11$}
\begin{center}
   \begin{tabular}{| c | c | c | c | c | c | }
     \hline
      & 1 & 2 & 3 & 4 & 5 \\ \hline
     5 & $A:5^2.2$ & $A:5^2.2$ & $A:5^2.2$ & $A:5^2.2$ & $A:2^5.1^2$ \\ \hline
     4 & Ir. & Ir. & Ir. & $A:5^2.1^2$ &\\
     &  &  &  & $B:6.3^2$ &\\
     \hline
     3 & Ir. & Ir. & $A:10.1^2$ & &  \\
          &  &  & $B:8.4$ & &  \\ \hline
     2 & Ir. & $A:5^2.1^2$ & & &  \\
      &  & $B:6^2$ & & &  \\
               \hline
    1 & $A:10.1^2$ & & & &   \\
                   \hline           
   \end{tabular}
   \end{center}
\vskip 0.5 cm
\centerline{Reducibility of $H$: $p=13$}
\begin{center}
   \begin{tabular}{| c | c | c | c | c | c | c |}
     \hline
      & 1 & 2 & 3 & 4 & 5 & 6\\ \hline
     6 & $A:6^2.2$ & $A:6.3^2.2$ & $A:4^2.2^3$ & $A:6.3^2.2$ & $A:6^2.2$ & $A:2^6.1^2$\\ 
     &  &  &  &  &  & $B:7^2$\\  \hline
     5 & Ir. & Ir. & Ir. & Ir. & $A:12.1^2$ &\\   \hline
     
     4 & $A:6.4^2$ & $A:6.4^2$ & Ir.&$A:3^4.1^2$ & & \\
       &  &  & & $B:7^2$ & &     \\ \hline
     3 & $A:8.6$ & Ir. & $A:4^3.1^2$ & & &  \\ \hline
     2 & Ir. & $A:6^2.1^2$ & & & &  \\
     & & $B:7^2$ & & & &  \\
               \hline
    1 & $A:12.1^2$ & & & & &   \\
                   \hline           
   \end{tabular}
   \end{center}
   \vskip 0.5 cm
 \par Our second conjecture is the following. Note that the irrationality measure for $\alpha$, given in both conjectures below, is directly derived applying Corollary 4.2, with $\ell=2k-j$.
 \newline {\bf{Conjecture 5.3. }}{\emph{Let $k\geq 1$ and $m\geq 2$ be integers. Let $p$ be a prime number such that $p=km+1$. Let $V=[p,k,k,\epsilon_A]$ and $H(V)$ be defined as above. Then the polynomial $H(V)$ has a factor $P$ in $\F_p[T][X]$ such that $P(\alpha)=0$, $\deg_X(P)=m$ and $\nu(\alpha)=m$. More precisely, we have
 $$P(X)=\sum_{0\leq i\leq m}(-1)^{ki}\binom{m}{i}T^{(1-(-1)^i)/2}X^{m-i}.$$
 }} 
\par Finally, we give a last conjecture allowing us to return to the quartic equation which has been the starting point of this investigation. 
  \newline {\bf{Conjecture 5.4. }}{\emph{Let $p$ be a prime number such that $p\equiv 1 \mod 3$ and $k=(p-1)/3$. Let $V=[p,k,k/2,\epsilon_A]$ and $H(V)$ be defined as above. Then the polynomial $H(V)$ has a factor $P$ in $\F_p[T][X]$ such that $P(\alpha)=0$, $\deg_X(P)=4$ and $\nu(\alpha)=8/3$. More precisely, with $j=k/2$, we have
  $$P(X)=X^4-w_{j+1}TX^3-(w_{j+1}/w_{j+2})X^2-(w_{j+1}/w_{j+2})^2/12.$$
  }} 
 \newline If $V=[13,4,2,-1]$ then we have $P=X^4+6TX^3+2X^2+4$. By elementary transformations on $P$ and $\alpha$ (see \cite[p.~1114-1115]{L3}), we return to the quartic equation $X^4+X^2-TX+1=0$, having a solution in $\F(13)$, with predictable continued fraction expansion, introduced by Mills and Robbins in \cite[p.~403]{MR}.
\par At last, we want to add that the interested reader can willingly obtain the code for the programs supporting these conjectures, by writing to the author.
\vskip 0.5 cm
\noindent {\bf{Aknowledgements. }} We would like to warmly thank our friend Nicolas Brisebarre for his steady and generous support during the preparation of this note. Let us also express our gratitude for the skillful advices on computer programming given by Bill Allombert.

\vskip 0.5 cm
\begin{tabular}{ll}Alain LASJAUNIAS\\Institut de Math\'ematiques de Bordeaux  CNRS-UMR 5251
\\Universit\'e de Bordeaux \\Talence 33405, France \\E-mail: Alain.Lasjaunias@math.u-bordeaux.fr\\\end{tabular}

\end{document}